\documentclass[11pt, righteqno]{article}
\usepackage{epsfig}
\usepackage{amsfonts}
\usepackage{amsmath}
\usepackage{amssymb}
\makeindex  \textwidth 130mm

\newtheorem{teo}{\sc Theorem}[section]
\newtheorem{proposition}{\sc Proposition}[section]

\newtheorem{lema}{\sc Lemma}[section]

\newtheorem{prf*}{\sc Proof}

\providecommand{\U}[1]{\protect\rule{.1in}{.1in}}

\providecommand{\U}[1]{\protect\rule{.1in}{.1in}}

\title{Quadratic Differential
Systems on $\mathbb{R}^{3}$ Having a Semisimple Derivation}
\author{I. Burdujan}

\begin{document}

\title{Quadratic Differential
Systems on $\mathbb{R}^{3}$ Having a non-Singular Semisimple
Derivation}
\author{I. Burdujan}

\date{}
\maketitle

\footnote{Corresponding author,

\noindent Ilie Burdujan, \emph{E-mail address}:\ $ilieburdujan@uaiasi.ro$, $ilieburdujan@yahoo.com$\\
address:\  3, Mihail Sadoveanu Street, 700490, Ia\c si, Romania}

\noindent\small{University of Agricultural Sciences and Veterinary
Medicine, Ia\c si, 700490, Romania}
\begin{abstract}\noindent The classification, up to a
center-affinity, of the homogeneous quadratic differential systems
defined on $\mathbb{R}^{3}$ that have at least a semisimple
nonsingular derivation, is achieved. It is proved that there exist
four affine-equivalence classes of such systems, only.
\end{abstract}

\section{Introduction}

Let us recall that the problem of classification up to an
affine-equivalence of homogeneous quadratic differential systems
is equivalent to the problem of classification up to an
isomorphism of commutative binary algebras defined on the same
ground space as the analyzed systems (e.g., see \cite{Ma1},
\cite{Wa1}, \cite{Bur1}).

Let $k$ be a perfect field and $A$ be a finite-dimensional vector
space over $k$. Recall that an endomorphism $T\in End\ A$ is said
to be semisimple if it is diagonalisable in an extension of $k$,
i.e. there exists a basis of $A$ consisting of eigenvectors of $T$
(see \cite{Lang}). In fact, $T$ is semisimple if and only if its
minimal polynomial is separable.

In the following, we restrict our interest to the case when $k$ is
$\mathbb{R}$ or $\mathbb{C}$.

In general, a derivation of a $k$-algebra is uniquely represented,
via its {\sc Jordan-Chevalley} decomposition, as the sum of a
semisimple derivation and a nilpotent derivation. Some classes of
real 3-dimensional commutative algebras having a derivation were
already classified, up to an isomorphism, in \cite{Bur2},
\cite{Bur3} and \cite{Bur4}; more exactly, there were respectively
obtained the isomorphism classes for algebras having either a
derivation with complex eigenvalues (\cite{Bur2}, \cite{Bur3}) or
a nilpotent derivation (\cite{Bur4}). Accordingly, there were
classified their corresponding homogeneous quadratic differential
systems up to a center-affine equivalence.

The aim of this paper is to classify, up to an isomorphism, the
real 3-dimensional commutative algebras having at least a
semisimple nonsingular derivation. The main result proved here is:
\emph{there exist four isomorphism classes of real 3-dimensional
commutative algebras having at least a semisimple nonsingular
derivation}. For each of them are exhibited the main properties
which allow to decide rapidly on the problem of their mutual
isomorphism relationships. Really, the existence of a semisimple
nonsingular derivation works like a strong constraint on the
integral curves of the associated HQDS; more exactly, all
nonsingular integral curves must be torsion free and, for a large
part of them the curvature tensor have to be zero as well. On the
other hand, from the algebraic point of view, the existence of a
semisimple nonsingular derivation compels the solvability of
analyzed algebras, yet the nilpotence of a large part of them.

In fact, we shall get the subalgebra lattices, the derivation
algebras and the group of automorphisms for each class of analyzed
algebras, because they are the most important invariants of any
binary algebra. Especially, we are interested in finding the set
$Ann\ A$ of annulator elements, the set $\mathcal{N}(A)$ of all
nilpotent elements and the set $\mathcal{I}(A)$ of all idempotent
elements of each analyzed algebra $A$. Further, the corresponding
homogeneous quadratic dynamical systems are classified up to a
center-affine equivalence. Finally, let us remark that subalgebra
lattice of $A(\cdot)$ allows to identify a natural partition
$\mathcal{P}_{A}$ of the ground space $A$ of algebra whose sets
will be named cells. In its turn, this partition has the property:
if an integral curve of the associated homogeneous quadratic
differential system meets a cell of partition $\mathcal{P}_{A}$
then it is entirely contained in that cell. It follows that the
set of all integral curves decomposes into a partition
subordinated to $\mathcal{P}_{A}$. Both these partitions are
invariant under the natural action of $Aut\ A$. They could be
useful tools in studying stability of singular solutions.

\section{Preliminaries}

Let $A$ be a finite dimensional vector space over a field $k$,
$k[A]$ be the coordinate ring of $A$,
$\mathcal{B}=(e_{1},e_{2},...,e_{n})$ a basis in $A$ and
$(x^{1},x^{2},...,x^{n})$ - the coordinate system assigned to
$\mathcal{B}$.

\vspace{2mm}Any autonomous differential system $(\mathcal{S})$ on
$A$ of the form
\begin{equation}\label{e1} \dot{x}^{i}=a_{jk}^{i}x^{j}x^{k}\ \ \ \ \ i,j,k\in \{1,2,...,n\}\end{equation}
with $a_{jk}^{i}=a_{kj}^{i}\in k$ and
$\dot{x}^{i}=\frac{dx^{i}}{dt}$ is called a \emph{homogeneous
quadratic differential system} (shortly, HQDS) on $A$. To any HQDS
(\ref{e1}) there is naturally associated a commutative binary
algebra $A(\cdot)$ defined by extending bi-linearly to $A\times A$
the mapping from $\mathcal{B}\times \mathcal{B}\rightarrow A$
defined by: $$(e_{j}, e_{k})\rightarrow e_{j}\cdot
e_{k}=a_{jk}^{i}e_{i}.$$ Conversely, the structure constants of
any finite dimensional commutative algebra could be used as
coefficients defining a system of type (\ref{e1}). This result
opens the way for an algebraic approach of HQDSs.

In fact, the structure properties of $A(\cdot)$ (i.e., properties
which are invariant under isomorphism action) enforce the
existence of some qualitative properties of ($\mathcal{S}$).
Conversely, the qualitative properties of $(\mathcal{S})$ induce
structure properties for $A(\cdot)$.

More exactly, there exists a correspondence (not necessarily a
map)
\begin{center} $\{$\emph{the properties affinely
invariant of any HQDE} $(\mathcal{S})$$\}$
$\leftrightsquigarrow$\end{center} \vspace{1mm}
$\leftrightsquigarrow$$\{$\emph{the properties invariant up to an
isomorphism of the associated algebra} $A(\cdot)$$\}$

\vspace{2mm} Unfortunately, this correspondence is not completely
discovered. Some of its components, revealed in \cite{KiSa1},
\cite{KiSa2}, \cite{Bur1} et al., are:

\vspace{2mm} 1. the set of nilpotents of $A(\cdot)$ is in a
bijective correspondence with the set of steady state (constant)
solutions of ($\mathcal{S}$),

\vspace{1mm} 2. the set of idempotents of $A(\cdot)$ is in a
bijective correspondence with the uniparametric sets of ray
solutions of ($\mathcal{S}$),

\vspace{1mm} 3. $A(\cdot)$ is a nilalgebra if and only if all
solutions of ($\mathcal{S}$) are polynomials,

\vspace{1mm} 4. $A^{2}=A\cdot A$ is a proper ideal of $A$ if and
only if ($\mathcal{S}$) has at least a \emph{linear first
integral},

\vspace{1mm} 5. if $A(\cdot)$ is a power-associative algebra then
all solutions of ($\mathcal{S}$) are rational functions,

\vspace{1mm} 6. if $A(\cdot)$ has a nonzero idempotent element
then zero solution of ($\mathcal{S}$) is not asymptotically
stable; consequently, if ($\mathcal{S}$) has a {\sc Liapunov}
function then $A$ has no idempotent element.

\vspace{3mm} Moreover, any proper ideal $I$ of $A$ decouples a
subsystem which is the HQDS associated with algebra $A/I$ such
that all nonsingular integral curves of $(\mathcal{S})$ lie on
cylinders. In the particular case when $Ann\ A$ (here $Ann\
A=\{v\in A\ |\ v\cdot w=0\ for\ all\ w\in A \}$) is a nonzero
ideal there exists at least a linear prime integral, such that all
nonsingular integral curves of $(\mathcal{S})$ are torsion-free.

\vspace{3mm} A main invariant object of any (commutative) algebra
$A$ is its subalgebra lattice $\mathcal{S}_{A}$. In fact, by using
some results due to {\sc Kaplan$\&$Yorke} \cite{KaYo1} and {\sc
Krasnov-Kononovich-Osharovich} \cite{KKOsh}, there will be proved
that, for any algebra $A$ analyzed in this paper, its subalgebra
lattice $\mathcal{S}(A)$ is fully determined by the set of its
1-dimensional subalgebras $\mathbb{R}u$ for $u\in
\mathcal{N}(A)\cup \mathcal{I}(A)$.

\vspace{3mm} The subalgebra lattice of $A(\cdot)$ allows to
identify a natural partition $\mathcal{P}(A)$ of $A$ consisting of
following nonempty parts of $A$:

\vspace{2mm}- the families of single sets covering each null
1-dimensional subalgebra,

- the connected components of each non-null 1-dimensional
subalgebra, delimited by $O$,

- the connected components of each 2-dimensional subalgebra
delimited by its subalgebra lattice of 1-dimensional subalgebras,

- the connected components of $A$ delimited by its subalgebra
lattice (i.e. the family of connected components of $A\setminus
\underset{{B\in \mathcal{S}_{A}}}{\bigcup}B$).

\vspace{3mm} Accordingly, there exists a natural partition of the
set of all integral curves for the corresponding HQDS.

\vspace{3mm} It was proved that there exists a 1-to-1 mapping
between affine-equivalent classes of HQDSs on $A$ and isomorphism
classes on commutative binary algebras on $A$. Consequently, the
problem of classification up to an affine-equivalence of HQDSs is
equivalent to the problem of finding the isomorphism classes of
commutative binary algebras. In order to solve this last problem,
it is suitable to divide the class of commutative algebras into
two disjoint parts:

\vspace{2mm}\ (i) the class of algebras having at least a
derivation,

(ii) the class of algebras having no derivation.

\vspace{2mm}For algebras belonging to former class there exist
more facilities for their studies than for the algebras in the
second class. Indeed, the existence of any derivation implies the
existence of a one-parametric group of automorphisms what assure
the existence of suitable bases where structure constants of
algebra have a simpler/simplest form, i.e. most of them are $0,\
1$ or other small integers. In the specific case under analysis in
this paper, the isomorphism classes have representatives which
have either $0$ or $1$ as structure constants in suitable bases.

\section{Algebras having a semisimple derivation}

\setcounter{equation}{0} \vspace{3mm}\noindent Let $A(\cdot)$ be
the real 3-dimensional (nontrivial) commutative algebra associated
with a HQDS in $A$. For convenience, suppose that the real vector
space $A$ was identified with $\mathbb{R}^{3}$ by means of any
chosen basis.

\vspace{3mm}\noindent Suppose that $\widetilde{D}$ is a nonzero
semisimple derivation of $A(\cdot)$. Then, algebra $A(\cdot)$ has
a derivation $D$ having the spectrum of the form $Spec\
D=(1,\lambda,\mu)$ where $1,\lambda,\mu$ are the eigenvalues,
distinct or not, of $D$. It follows that there exists a basis
$\mathcal{B}=(e_{1}, e_{2}, e_{3})$ of $A$ such that
$$D(e_{1})=e_{1},\ \ \ D(e_{2})=\lambda e_{2},\ \ \ D(e_{3})=\mu e_{3}.$$
Since in basis $\mathcal{B}'=(f_{1}=e_{1}, f_{2}=e_{3},
f_{3}=e_{2})$ we get $D(f_{1})=f_{1},\ D(f_{2}) =\mu f_{2},\
D(f_{3})=\lambda f_{3}$, it results it is enough to analyze only
the case $Spec\ D=(1,\lambda,\mu)$ with $\lambda\leq \mu$. Since
the derivation is nonsingular it results necessarily $\lambda\cdot
\mu\neq 0$. Of course, $A$ decomposes into a direct vector sum of
eigenspaces of derivation $D$.

In order to express analytically the existence of a derivation $D$
for algebra $A(\cdot)$, we consider the structure constants of $A$
in basis $\mathcal{B}$ usually defined by
\begin{equation}\label{e21}e_{i}\cdot e_{j}=a_{ij}^{k}e_{k}.\end{equation} For convenience, we shall denote
\begin{equation}\label{e22}\begin{array}{llllll}a_{11}^{1}=a&\ \ a_{11}^{2}=b&\
\ a_{11}^{3}=c&\ \ a_{12}^{1}=k&\ \ a_{12}^{2}=m&\ \ a_{12}^{3}=n\\

a_{22}^{1}=d&\ \ a_{22}^{2}=e&\
\ a_{22}^{3}=f&\ \ a_{13}^{1}=p&\ \ a_{13}^{2}=q&\ \ a_{13}^{3}=r\\

a_{33}^{1}=g&\ \ a_{33}^{2}=h&\ \ a_{33}^{3}=j&\ \ a_{23}^{1}=s&\
\ a_{23}^{2}=t&\ \ a_{23}^{3}=v.
\end{array}\end{equation}

\noindent Then $D$ is a derivation for $A$ if and only if the next
conditions are fulfilled:

\begin{equation}\label{e23}\left\{\begin{array}{l}
a=e=j=k=m=p=r=t=v=0\\
(\lambda-2)b=0\\
(\mu-2)c=0\\
(1-2\lambda)d=0\\
(\mu-2\lambda)f=0\\
(1-2\mu)g=0\\
(\lambda-2\mu)h=0\\
(\mu-1-\lambda)n=0\\
(\lambda-1-\mu)q=0\\
(1-\lambda-\mu)s=0.
\end{array}\right.\end{equation}
These equations impose to take into account of the natural
decomposition of the plane $\mathbb{R}^{2}$ of variables
$(\lambda,\ \mu)$ defined by means of the lines:
\begin{equation}\label{e24}\begin{array}{l}\lambda=2,\ \ \mu=2,\ \ 1-2\lambda=0,\ \ 1-2\mu=0,\ \ \mu-2\lambda=0,\\  \lambda-2\mu=0,
\mu-1-\lambda=0,\ \ \lambda-1-\mu=0,\ \
\lambda+\mu-1=0.\end{array}\end{equation}

The next result holds true.
\begin{lema} The set $\mathbb{R}^{2}\setminus {\{\lambda=0\}\cup \{\mu=
0\}}$ decomposes into the following disjoint subsets:

\vspace{3mm}i) on the line $\lambda=\frac{1}{2}$ :
$\{(\frac{1}{2},\ -\frac{1}{2})\},\ \{(\frac{1}{2},\
\frac{1}{4})\},\ \{(\frac{1}{2},\ \frac{1}{2})\},\
\{(\frac{1}{2},\ 1)\},\ \{(\frac{1}{2},\ \frac{3}{2})\},\
\{(\frac{1}{2},\ 2)\}$ and $\{(\frac{1}{2},\ \mu)\ |\ \mu\notin \{
-\frac{1}{2},\ 0,\ \frac{1}{4},\ \frac{1}{2},\ 1,\ \frac{3}{2},\
2\}\}$

ii) on the line $\lambda=2$ : $\{(2,-1)\},\ \{(2,1)\},\
\{(2,2)\},\ \{(2,3)\},\ \{(2,4)\},\  \{(2,\frac{1}{2})\}$, and
$\{(2,\ \mu)\ |\ \mu\notin\{-1,0,\frac{1}{2},1,2,3\}\}$

iii) on the line $\mu=\frac{1}{2}$ : $\{(-\frac{1}{2},\
\frac{1}{2})\},\ \{(\frac{1}{4},\ \frac{1}{2})\},\
\{(\frac{1}{2},\ \frac{1}{2})\},\ \{(1,\ \frac{1}{2})\},\
\{(\frac{3}{2},\ \frac{1}{2})\},\ \{(2,\ \frac{1}{2})\}$ and
$\{(\lambda,\ \frac{1}{2})\ |\ \lambda\notin \{-\frac{1}{2},\ 0, \
\frac{1}{4},\ \frac{1}{2},\ 1,\ \frac{3}{2},\ 2\}\}$

iv) on the line $\mu=2$ : $\{(-1,2)\},\ \{(1,2)\},\ \{(2,2)\},\
\{(3,2)\},\ \{(4,2)\},\  \{(\frac{1}{2},2)\}$, and $\{(\lambda,\
2)\ |\ \lambda\notin\{-1,0,\frac{1}{2},1,2,3,4\}\}$

v) on the line $\lambda=2\mu$ : $\{(-2,-1)\},\ \{(\frac{1}{2},
\frac{1}{4})\},\ \{(\frac{2}{3},\ \frac{1}{3})\},\
\{(1,\frac{1}{2})\},\ \{(2,1)\},\ \{(4,2)\},$ and $\{(2\mu,\ \mu)\
|\ \mu\notin \{-2,\ 0,\ \frac{1}{4}, \ \frac{1}{3}, \
\frac{1}{2},\ 1,\ 2 \}\}$

vi) on the line $\mu=\lambda-1$ : $\{(-1,-2)\},\
\{(\frac{1}{2},-\frac{1}{2})\},\ \{(\frac{3}{2},\frac{1}{2})\},\
\{(2,1)\},\ \{(3,2)\},$ and $\{(\lambda,\ \lambda-1)\ |\
\lambda\notin\{-1,\ 0,\ \frac{1}{2},\ \frac{3}{2},\ 2,\ 3,\}\}$

vii) on the line $\mu=\lambda+1$ : $\{(-2,-1)\},\
\{(-\frac{1}{2},\frac{1}{2})\},\ \{(\frac{1}{2},\frac{3}{2})\},\
\{(1,2)\},\ \{(2,3)\},$ and $\{(\lambda,\ \lambda+1)\ |\
\lambda\notin\{-2,\ -\frac{1}{2},\ 0,\ \frac{1}{2},\ 1,\ 2\}\}$

viii)  on the line $\mu=2\lambda$ : $\{(-1,-2)\},\
\{(\frac{1}{4},\frac{1}{2})\},\ \{(\frac{1}{2},1)\},\ \{(1,2)\},\
\{(2,4)\},$ and $\{(\lambda,\ 2\lambda)\ |\ \lambda\notin\{-1,\
0,\ \frac{1}{4},\ \frac{1}{2},\ 1,\ 2\}\}$

ix) on the line $\mu=1-\lambda$ : $\{(-1,2)\},\
\{(\frac{1}{3},\frac{2}{3})\},\ \{(\frac{1}{2},\frac{1}{2})\},\
\{(\frac{2}{3},\frac{1}{3})\},\ \{(2,-1)\},$ and $\{(\lambda,\
1-\lambda)\ |\ \lambda\notin\{-1,\ 0,\ \frac{1}{3},\ \frac{1}{2},\
\frac{2}{3},\ 2\}\}$

x) $\mathbb{R}^{2}\setminus {\{\lambda=0\}\cup \{\mu= 0\}}$ less
the lines listed in (\ref{e24}).
\end{lema}

Accordingly, the following results are proved.

\begin{proposition}\label{p1} Any real 3-dimensional commutative algebra $A(\cdot)$
having at least a semisimple derivation with kernel $\{0\}$ has a
derivation $D$ such that its spectrum $Spec\ D=(1,\lambda,\mu)$
has as $(\lambda, \mu)$ one of the pairs listed in Lemma 3.1,
i)-ix).
\end{proposition}

\begin{proposition}\label{p2}  Let $A(\cdot)$ be any real 3-dimensional commutative
algebra having at least a semisimple derivation with kernel
$\{0\}$. Then $Der\ A$ contains at least a derivation having the
spectrum in one (and only one) of the next nine families

\vspace{3mm} 1) $\{(1,\frac{1}{2},\ -\frac{1}{2}),\
(1,-\frac{1}{2},\ \frac{1}{2}),\ (1,-1,2),\ (1,2,-1), (1,-2,-1),\
(1,-1,-2)\}$,

2) $\{(1,\frac{1}{2},\ \frac{1}{4}),\ (1,\frac{1}{4},\
\frac{1}{2}),\ (1,2,4),\ (1,4,2),\ (1,\frac{1}{2}, 2),\
(1,2,\frac{1}{2}) \}$,

3) $\{(1,\frac{1}{2},\ \frac{1}{2}),\ (1,2,1),\ (1,1,2)\}$,

4) $\{(1,\frac{1}{2},1),\ (1,1,\frac{1}{2}),\ (1,2,2)\}$,

5) $\{(1,\frac{1}{2},\frac{3}{2}),\ (1,\frac{3}{2},\frac{1}{2}),\
(1,2,3),\ (1,3,2),\ (1,\frac{2}{3},\frac{1}{3}),\
(1,\frac{1}{3},\frac{2}{3})\}$,

6) $\{(1,\frac{1}{2},\mu)\ |\ \mu\in \mathbb{R}\setminus \{
-\frac{1}{2},\ 0,\ \frac{1}{4},\ \frac{1}{2},\ 1,\ \frac{3}{2},\
2\}\}\ \cup \{(1,2,\ \mu)\ |\ \mu\in \mathbb{R}\setminus
\{-1,0,\frac{1}{2},1,2,3,4\}\}\cup \{(1,\lambda,\ \frac{1}{2})\ |\
\lambda\in \mathbb{R}\setminus \{-\frac{1}{2},\ 0, \ \frac{1}{4},\
\frac{1}{2},\ 1,\ \frac{3}{2},\ 2\}\}\cup  \{(1,\lambda,\ 2)\ |\
\lambda\in \mathbb{R}\setminus \{-1,0,\frac{1}{2},1,2,3,4\}\}$,

7) $\{(1, 2\mu,\ \mu)\ |\ \mu\in \mathbb{R}\setminus\{-1,\ 0,\
\frac{1}{4}, \ \frac{1}{3}, \ \frac{1}{2},\ 1,\ 2 \}\}\cup
\{(1,\lambda,\ 2\lambda)\ |\ \lambda\in \mathbb{R}\setminus \{-1,\
0, \ \frac{1}{4},\ \frac{1}{2},\ 1,\ 2\}\}$,

8) $\{(1,\lambda,\ \lambda-1)\ |\ \lambda\in
\mathbb{R}\setminus\{-1,\ 0, \ \frac{1}{2},\ \frac{3}{2},\ 2,\
3\}\}\cup \{(1,\lambda,\ \lambda+1)\ |\ \lambda\in
\mathbb{R}\setminus\{-2,\ -\frac{1}{2},\ 0,\ \frac{1}{2},\ 1,\
2\}\}$,

9) $\{(1,\lambda,\ 1-\lambda)\ |\ \lambda\in\mathbb{R}\setminus
\{-1,\ 0,\ \frac{1}{3},\ \frac{1}{2},\ \frac{2}{3},\ 2\}\}.$
\end{proposition}

Let us recall that whenever $D\in Der\ A$ has $Spec\
D=(\alpha,\beta,\gamma)$ then $\tilde{D}=kD$ is a derivation with
$Spec\ \tilde{D}=(k\alpha,\ k\beta,\ k\gamma)$. This remark and
Proposition \ref{p2} allow to prove the next result.

\begin{proposition}\label{p3} In order to classify the real 3-dimensional commutative algebras having a semisimple nonsingular derivation
it is enough to restrict the study on the algebras having a
derivation with spectrum of one of the following types:
$$\begin{array}{c}(1,-1,2),\ (1,1,2),\ (1,2,2),\ (1,2,3),\ (1,2,4)\\
(1,2,\ \lambda)\ with \ \lambda\in \mathbb{R}\setminus
\{-1,0,\frac{1}{2},1,2,3,4\},\\ (1,2,\ 2\lambda)\ with \
\lambda\in \mathbb{R}\setminus \{-\frac{1}{2},\ 0,\ \frac{1}{4},\
\frac{1}{2},\ 1,\ \frac{3}{2},\
2\}\},\\
(1,\lambda,\ 2\lambda)\ with \ \lambda\in \mathbb{R}\setminus
\{-1,\ \frac{1}{4},\ \frac{1}{2},\ 1,\ 2\},\\ (1,\lambda,\
\lambda+1)\ with \ \lambda\in
\mathbb{R}\setminus\{-2,\ -\frac{1}{2},\ \frac{1}{2},\ 1,\ 2\},\\
\{(1,\lambda,\ 1-\lambda)\ |\ \lambda\in\mathbb{R}\setminus \{-1,\
\frac{1}{3},\ \frac{1}{2},\ \frac{2}{3},\ 2\}.\end{array}$$
\end{proposition}

\vspace{3mm} In the following we will identify all these classes
of algebras.

\vspace{3mm} 1)\ \ \textbf{Case} $Spec\ D=(1,-1,2)$

\vspace{3mm}\noindent There exists a basis
$\mathcal{B}=(e_{1},e_{2},e_{3})$ such that the multiplication
table of algebra $A(\cdot)$ has the form:

 $$\begin{array}{llll}
  \textbf{Table T}\hspace{6mm} & \hspace{5mm} e_{1}^{2}=pe_{3}  &\hspace{5mm} e_{2}^{2}=0&\hspace{5mm} e_{3}^{2}=0 \\
   &\hspace{5mm}  e_{1}e_{2}=0  &\hspace{5mm}  e_{1}e_{3}=0 &\hspace{5mm}  e_{2}e_{3}=qe_{1}
\end{array}$$
with $p,q\in\mathbb{R}$.

\vspace{3mm} \textbf{Subcase }$pq\neq 0$

\vspace{3mm}In the case when $pq\neq 0$ any such algebra is
isomorphic to algebra:
$$\begin{array}{llll}
  \textbf{Table T1}\hspace{6mm} & \hspace{5mm} e_{1}^{2}=e_{3}  &\hspace{5mm} e_{2}^{2}=0&\hspace{5mm} e_{3}^{2}=0 \\
   &\hspace{5mm}  e_{1}e_{2}=0&\hspace{5mm}  e_{1}e_{3}=0  &\hspace{5mm}  e_{2}e_{3}=e_{1}
\end{array}$$
(indeed, it is enough to use the basis
$(\frac{1}{p}e_{1},\frac{1}{q}e_{2},\frac{1}{p}e_{3})$).

\vspace{3mm} \emph{Properties of algebra $A_{1}\ (with\ table\
T1)$:}

\vspace{2mm} $\bullet$ $Ann\ A=\{0\},\
\mathcal{N}(A)=\mathbb{R}e_{2}\cup \mathbb{R}e_{3},\
\mathcal{I}(A)=\emptyset$,

$\bullet$ 1-dimensional subalgebras: $\mathbb{R}e_{2},\
\mathbb{R}e_{3}$,

$\bullet$ 2-dimensional subalgebras: $Span_{\mathbb{R}}\{e_{1},
e_{3}\}$,

$\bullet$ ideals: $ Span_{\mathbb{R}}\{e_{1}, e_{3}\}$,

$\bullet$ $A^{2}=Span_{\mathbb{R}}\{e_{1}, e_{3}\}$,

$\bullet$ $A$ is solvable; $A= \mathbb{R}e_{2}\oplus
Span_{\mathbb{R}}\{e_{1}, e_{3}\}$ - vector direct sum of a null
subalgebra and a nilpotent (maximal) ideal, i.e. this is a {\sc
Weddeburn-Artin} decomposition for $A$,

$\bullet$ $Der\ A= \mathbb{R}D$,

$\bullet$ $Aut\ A= \left\{
\left[\begin{array}{lll}x&0&0\\0&x^{-1}&0\\0&0&x^{2}\end{array}\right]\
|\ x\in \mathbb{R}^{\ast}\right\},$

$\bullet$ the partition $\mathcal{P}_{A}$ of $\mathbb{R}^{3}$,
defined by the subalgebra lattice of $A$, consists of:

\vspace{1mm} \hspace{5mm} $\diamond$  the singletons covering the
axes $Ox^{2}$ and $Ox^{3}$,

\hspace{5mm} $\diamond$  the half-planes delimited by axis
$Ox^{3}$ on $x^{1}Ox^{3}$,

\hspace{5mm} $\diamond$ the half-spaces $x^{2}>0$ and $x^{2}<0$
(delimited by plane $x^{1}Ox^{3}$) less the points of axis
$Ox^{2}$,

$\bullet$ the partition $\mathcal{P}_{A}$ of $A$ induces a
partition of the set of integral curves of the associated
homogeneous quadratic differential system (HQDS) consisting of:

\vspace{1mm} \hspace{5mm} $\diamond$  the singletons consisting of
singular solutions that cover the axes $Ox^{2}$ and $Ox^{3}$,

\hspace{5mm} $\diamond$  the integral curves contained in each
half-plane delimited by axis $Ox^{3}$ on $x^{1}Ox^{3}$,

\hspace{5mm} $\diamond$ the integral curves contained in each
half-space delimited by plane $x^{1}Ox^{3}$ which do not meet the
axis $Ox^{2}$,

The presence of ideal $Span_{\mathbb{R}}\{e_{1}, e_{3}\}$ assures
the existence of a decoupled subsystem defined on the null algebra
$\mathbb{R}e_{2}\cong A/Span_{\mathbb{R}}\{e_{1}, e_{3}\}$.

Since $A$ is solvable, the associated HQDS has a linear prime
integral and, consequently, each nonsingular integral curve is
torsion-free.

Note that algebra $A$ in neither associative nor
power-associative.

\vspace{3mm} \textbf{Subcase }$pq = 0$

\vspace{3mm}Of course, the case when $p=q=0$ gives the null
algebra which is not of any interest. It remains to analyze the
cases
$$(i)\ \ p\neq 0,\ \ q=0,\ \ \ \ (ii)\ \ p=0,\ \ q\neq 0.$$

\vspace{3mm} \textbf{(i)} In this case any such algebra is
isomorphic to algebra:

$$\begin{array}{llll}
  \textbf{Table T'2}\hspace{6mm} & \hspace{5mm} e_{1}^{2}=e_{3}  &\hspace{5mm} e_{2}^{2}=0&\hspace{5mm} e_{3}^{2}=0 \\
   &\hspace{5mm}  e_{1}e_{2}=0&\hspace{5mm}  e_{1}e_{3}=0  &\hspace{5mm}
   e_{2}e_{3}=0
\end{array}$$
Changing the basis $(e_{1},e_{2},e_{3})$ in basis
$(e_{2},e_{3},e_{1})$ the multiplication table of algebra
$A(\cdot)$ gets the form:
$$\begin{array}{llll}
  \textbf{Table T2}\hspace{6mm} & \hspace{5mm} e_{1}^{2}=0  &\hspace{5mm} e_{2}^{2}=0&\hspace{5mm} e_{3}^{2}=e_{2} \\
   &\hspace{5mm}  e_{1}e_{2}=0&\hspace{5mm}  e_{1}e_{3}=0  &\hspace{5mm}
   e_{2}e_{3}=0
\end{array}$$

\vspace{3mm} \emph{Properties of algebra $A_{2}\ (with\ table\
T2)$:}

\vspace{2mm}$\bullet$ $Ann\ A=Span_{\mathbb{R}}\{e_{1},e_{2}\},\
\mathcal{N}(A)=Span_{\mathbb{R}}\{e_{1},e_{2}\},\
\mathcal{I}(A)=\emptyset$,

$\bullet$ 1-dimensional subalgebras: $\mathbb{R}u$ for each $u\in
\mathcal{N}(A)$,

$\bullet$ 2-dimensional subalgebras: $Span_{R}\{e_{2},
pe_{1}+qe_{3}\}$ ($p^{2}+q^{2}\neq 0$),

$\bullet$ ideals: $\mathbb{R}u$ for $u\in
Span_{\mathbb{R}}\{e_{1},e_{2}\}$, $Span_{R}\{e_{2},
pe_{1}+qe_{3}\}$ ($p^{2}+q^{2}\neq 0$),

$\bullet$ $A^{2}=\mathbb{R}e_{2}$,

$\bullet$ $A$ is a nilpotent associative algebra,

$\bullet$ $Der\ A=\left\{
\left[\begin{array}{rrr}x&0&u\\y&2z&v\\0&0&z\end{array}\right]\ |\
x,y,z,u,v\in \mathbb{R}\right\}$,

$\bullet$ $Aut\ A= \left\{
\left[\begin{array}{rrr}x&0&u\\y&z^{2}&v\\0&0&z\end{array}\right]\
|\ x,y,z,u,v\in \mathbb{R},\ xz\neq 0\right\}$.

$\bullet$ the partition $\mathcal{P}_{A}$ of $\mathbb{R}^{3}$,
defined by the lattice of subalgebras of $A$, consists of:

\vspace{1mm}\hspace{5mm} $\diamond$  the singletons covering the
plane $x^{1}Ox^{2}$,

\hspace{5mm} $\diamond$  the half-planes delimited by axis
$Ox^{2}$ on each plane containing $Ox^{2}$ less $x^{1}Ox^{2}$,

$\bullet$ the partition $\mathcal{P}_{A}$ of $A$ induces a
partition of the set of integral curves of the associated HQDS
consisting of:

\vspace{1mm}\hspace{5mm} $\diamond$  the singletons consisting of
singular solutions covering the plane $x^{1}Ox^{2}$,

\hspace{5mm} $\diamond$  the integral curves contained in each
half-plane delimited by axis $Ox^{2}$ on any plane passing through
$Ox^{2}$ less $x^{1}Ox^{2}$.

As each nonsingular integral curve lies on a 2-dimensional
subalgebra, it results its torsion is zero. The presence of a
large family of ideals implies the existence of corresponding
family of prime integrals and, consequently, each nonsingular
integral curve has also a zero curvature tensor. More exactly,
each nonsingular integral curve lies on a parallel to $Ox^{2}$.

\vspace{3mm} \textbf{(ii)} In tis case, by using the changes of
bases $(e_{1},\frac{1}{q}e_{2},e_{3})$ and $(e_{3},e_{2},e_{1})$,
it follows that any analyzed algebra is isomorphic to algebra:

$$\begin{array}{llll}
  \textbf{Table T3}\hspace{6mm} & \hspace{5mm} e_{1}^{2}=0  &\hspace{5mm} e_{2}^{2}=0&\hspace{5mm} e_{3}^{2}=0 \\
   &\hspace{5mm}  e_{1}e_{2}=e_{3}&\hspace{5mm}  e_{1}e_{3}=0  &\hspace{5mm}
   e_{2}e_{3}=0
\end{array}$$

\vspace{3mm} \emph{Properties of algebra $A_{3}\ (with\ table\
T3)$:}

\vspace{2mm}$\bullet$ $Ann\ A=\mathbb{R}e_{3},\
\mathcal{N}(A)=Span_{\mathbb{R}}\{e_{1},e_{3}\}\cup
Span_{\mathbb{R}}\{e_{2},e_{3}\},\ \mathcal{I}(A)=\emptyset$,

$\bullet$ 1-dimensional subalgebras: $\mathbb{R}u$ for each $u\in
\mathcal{N}(A)$,

$\bullet$ 2-dimensional subalgebras: $Span_{R}\{e_{3},
pe_{1}+qe_{2}\}$ ($p^{2}+q^{2}\neq 0$),

$\bullet$ ideals: $\mathbb{R}e_{3},\ Span_{R}\{e_{3},
pe_{1}+qe_{2}\}$ ($p^{2}+q^{2}\neq 0$),

$\bullet$ $A^{2}=\mathbb{R}e_{3}$,

$\bullet$ $A$ is a nilpotent associative algebra,

$\bullet$ $Der\ A=\left[\begin{array}{rrc}x&0&0\\
0&z&0\\y&v&x+z\end{array}\right]$ with $x,y,z,v\in \mathbb{R}$,

$\bullet$ $Aut\ A= H\cup h\cdot H$ where
$$H= \left\{
\left[\begin{array}{rrc}x&0&0\\0&z&0\\y&v&xz\end{array}\right]\ |\
x,y,z,v\in \mathbb{R},\ xz\neq 0\right\}\ \ \text{and}\ \ h=
\left[\begin{array}{rrr}0&1&0\\ 1&0&0\\0&0&1\end{array}\right]$$
(in fact, $H$ is a normal divisor for $Aut\ A$, $h^{2}=Id$ and
$Aut\ A/H \cong \mathbb{Z}_{2}$),

$\bullet$ the partition $\mathcal{P}_{A}$ of $\mathbb{R}^{3}$,
defined by the lattice of subalgebras of $A$, consists of:

\vspace{1mm}\hspace{5mm} $\diamond$  the singletons covering the
planes $x^{1}Ox^{3}$ and $x^{2}Ox^{3}$,

\hspace{5mm} $\diamond$  the half-planes delimited by axis
$Ox^{3}$ on each plane containing $Ox^{3}$ less the planes
$x^{1}Ox^{3}$ and $x^{2}Ox^{3}$,

$\bullet$ the partition $\mathcal{P}_{A}$ of $A$ induces a
partition of the set of integral curves of the associated HQDS
consisting of:

\vspace{1mm}\hspace{5mm} $\diamond$  the singletons consisting of
singular solutions covering the planes $x^{1}Ox^{3}$ and
$x^{2}Ox^{3}$,

\hspace{5mm} $\diamond$ the integral curves lying on the
half-planes delimited by axis $Ox^{3}$ of each plane containing
$Ox^{3}$ less the planes $x^{1}Ox^{3}$ and $x^{2}Ox^{3}$.

Each nonsingular integral curve lies on a subalgebra so that it is
torsion-free. There exists a large family of linear prime
integrals in accordance with the family of ideals and,
consequently, each nonsingular integral curve has a zero curvature
tensor. More exactly, each nonsingular integral curve lies on a
parallel to $Ox^{3}$.

\vspace{3mm}2) \textbf{Case} $Spec\ D=(1,1,2)$

\vspace{3mm}\noindent There exists a basis
$\mathcal{B}=(e_{1},e_{2},e_{3})$ such that the multiplication
table of algebra $A(\cdot)$ has the form:

 $$\begin{array}{llll}
  \textbf{Table T}\hspace{6mm} & \hspace{5mm} e_{1}^{2}=pe_{3}  &\hspace{5mm} e_{2}^{2}=qe_{3}&\hspace{5mm} e_{3}^{2}=0 \\
   &\hspace{5mm}  e_{1}e_{2}=re_{3}  &\hspace{5mm}  e_{1}e_{3}=0 &\hspace{5mm}
   e_{2}e_{3}=0
\end{array}$$
with $p,q,r\in \mathbb{R}$.

\vspace{3mm}\textbf{Subcase  $pqr\neq 0$}

\vspace{3mm}\noindent Each such algebra is isomorphic to algebra:
$$\begin{array}{llll}
  \textbf{Table T'}\hspace{6mm} & \hspace{5mm} e_{1}^{2}=e_{3}  &\hspace{5mm} e_{2}^{2}=\lambda e_{3} &\hspace{5mm} e_{3}^{2}=0 \\
   &\hspace{5mm}  e_{1}e_{2}=e_{3} &\hspace{5mm}  e_{1}e_{3}=0  &\hspace{5mm}
   e_{2}e_{3}=0
\end{array}$$
with $\lambda\neq 0$ (indeed, it is enough to use the basis
$(e_{1},\frac{p}{r}e_{2},pe_{3})$). We have to consider the
subcases:
$$(i^{\circ})\  \lambda=1,\ \ \ \ (ii^{\circ})\  \lambda<1\ (\lambda\neq 0),\ \ \ \ (iii^{\circ})\  \lambda>1$$

\vspace{3mm}$(i^{\circ})$\  \textbf{Subcase}\ $\lambda=1$

\vspace{3mm}\noindent There exists a basis
$\mathcal{B}=(e_{1},e_{2},e_{3})$ such that the multiplication
table of algebra $A(\cdot)$ gets the form \textbf{T2}, so that
this algebra is isomorphic to algebra $A_{2}$.

\vspace{3mm}$(ii^{\circ})$ \textbf{Subcase} $\lambda<1,\ \lambda\
\neq 0$

\vspace{3mm}\noindent Let $s_{1}\neq s_{2}$ be the distinct
solutions of equation $\lambda s^{2}+2s+1=0$. Then, in basis
$(e_{1}+s_{1}e_{2},\ e_{1}+s_{2}e_{2},\
2\frac{\lambda-1}{\lambda}e_{3})$ where $\lambda
s_{i}^{2}+2s_{i}+1=0$ (for $i=1,2$) the multiplication table of
algebra $A(\cdot)$ gets the form \textbf{T3}, so that this algebra
is isomorphic to algebra $A_{3}$.

\vspace{3mm}$(iii^{\circ})$ \textbf{Subcase} $\lambda>1$

\vspace{3mm}\noindent In basis $(\sqrt{\lambda-1}e_{1},\
e_{1}-e_{2},\ (\lambda-1)e_{3})$ the multiplication table of
algebra $A(\cdot)$ gets the form:

 $$\begin{array}{llll}
  \textbf{Table T4}\hspace{6mm} & \hspace{5mm} e_{1}^{2}=e_{3}  &\hspace{5mm} e_{2}^{2}=e_{3}&\hspace{5mm} e_{3}^{2}=0 \\
   &\hspace{5mm}  e_{1}e_{2}=0  &\hspace{5mm}  e_{1}e_{3}=0 &\hspace{5mm}
   e_{2}e_{3}=0
\end{array}$$

\vspace{3mm} \emph{Properties of algebra $A_{4}\ (with\ table\
T4)$:}

\vspace{2mm}$\bullet$ $Ann\ A=\mathbb{R}e_{3},\
\mathcal{N}(A)=\mathbb{R}e_{3},\ \mathcal{I}(A)=\emptyset$,

$\bullet$ 1-dimensional subalgebras: $\mathbb{R}e_{3}$,

$\bullet$ 2-dimensional subalgebras: $Span_{R}\{e_{3},
pe_{1}+qe_{2}\}$ ($p^{2}+q^{2}\neq 0$),

$\bullet$ ideals: $\mathbb{R}e_{3},\ Span_{R}\{e_{3},
pe_{1}+qe_{2}\}$ ($p^{2}+q^{2}\neq 0$),

$\bullet$ $A^{2}=\mathbb{R}e_{3}$,

$\bullet$ $A$ is a nilpotent associative algebra,

$\bullet$ $Der\ A=\left[\begin{array}{rrr}x&-y&0\\
y&x&0\\z&v&2x\end{array}\right]$ with $x,y,z,v\in \mathbb{R}$,

$\bullet$ $Aut\ A=  \left\{ \left[\begin{array}{ccc}\rho cos\
\theta&-\rho sin\ \theta&0\\ \rho sin\ \theta&\rho cos\
\theta&0\\z&v&\rho^{2}\end{array}\right]\ |\ \rho ,z,v\in
\mathbb{R},\ \rho> 0,\ \theta\in[0,2\pi)\right\}$,

$\bullet$ the partition of $\mathbb{R}^{3}$, defined by the
lattice of subalgebras of $A$, consists of:

\vspace{1mm}\hspace{5mm} $\diamond$  the singletons covering axis
$Ox^{3}$,

\hspace{5mm} $\diamond$  the half-planes delimited by axis
$Ox^{3}$ on each plane containing $Ox^{3}$,

$\bullet$ the partition $\mathcal{P}_{A}$ of $A$ induces a
partition of the set of integral curves of the associated HQDS
consisting of:

\vspace{1mm}\hspace{5mm} $\diamond$  the singletons consisting of
singular solutions covering the axis $Ox^{3}$,

\hspace{5mm} $\diamond$ the integral curves lying on the
half-planes delimited by axis $Ox^{3}$ of each plane containing
$Ox^{3}$.

Each nonsingular integral curve lies on a subalgebra such that it
is torsion-free. There exists a large family of linear prime
integrals in accordance with the family of ideals and,
consequently, each nonsingular integral curve has a zero curvature
tensor. More exactly, each nonsingular integral curve lies on a
parallel to $Ox^{3}$.

\vspace{3mm}\textbf{Subcase  $pqr= 0$}

\vspace{3mm}\noindent In this case we have to analyze the
following six possibilities:

$$\begin{array}{llllll}(1)&\ p= 0,\ q\neq 0,\ r=0,& (2) & \ p= 0,\ q= 0,\
r\neq 0,& (3)&\ p= 0,\ q\neq 0,\ r\neq 0,\\ (4) & \ p\neq 0,\ q=
0,\ r= 0,& (5)&\ p\neq 0,\ q= 0,\ r\neq 0,& (6) & \ p\neq 0,\
q\neq 0,\ r= 0.
\end{array}$$

\vspace{3mm}\noindent \textbf{(1)}  This algebra is isomorphic to
algebra $A_{2}$.

\vspace{3mm}\noindent \textbf{(2)}  This algebra is isomorphic to
algebra $A_{3}$.

\vspace{3mm}\noindent \textbf{(3)}  This algebra is isomorphic to
algebra $A_{3}$.

\vspace{3mm}\noindent \textbf{(4)}  This algebra is isomorphic to
algebra $A_{2}$.

\vspace{3mm}\noindent \textbf{(5)} By using the changes of bases
$(\frac{r}{p}e_{1},e_{2},\frac{r^{2}}{p}e_{3})$ followed by
$(2e_{1}-e_{2},e_{2},2e_{3})$ it results that this algebra is
isomorphic to algebra $A_{3}$.

\vspace{3mm}\noindent \textbf{(6)} If $pq>0$ this algebra is
isomorphic to algebra $A_{4}$. In case $pq<0$ there exists a basis
such that the multiplication table of algebra $A(\cdot)$ gets the
form:

 $$\begin{array}{llll}
  \textbf{Table T5}\hspace{6mm} & \hspace{5mm} e_{1}^{2}=e_{3}  &\hspace{5mm} e_{2}^{2}=-e_{3}&\hspace{5mm} e_{3}^{2}=0 \\
   &\hspace{5mm}  e_{1}e_{2}=0  &\hspace{5mm}  e_{1}e_{3}=0 &\hspace{5mm}
   e_{2}e_{3}=0
\end{array}$$
This algebra has $\mathcal{N}(A)=\{x(e_{1}\pm e_{2})+ze_{3}\ |\
x,z\in \mathbb{R}\}$. Then the basis $(\frac{1}{2}(e_{1}+
e_{2}),\frac{1}{2}(e_{1}- e_{2}),\frac{1}{2}e_{3})$ assures us
that $A(\cdot)$ is isomorphic with an algebra of type $A_{3}$.

\vspace{5mm} \textbf{Case }$Spec\ D=(1,2,2)$

\vspace{3mm}\noindent Looking at the set of derivations for
algebra T2 we see that it contains a derivation $D$ with $Spec\
D=(1,2,2)$, namely, the derivation corresponding to
$$\beta=2,\ \ \rho=1,\ \ \gamma=\lambda=\mu=0;$$ accordingly, the
class of algebras having a derivation with spectrum $(1,2,2)$ is
contained into the class of algebras of type \textbf{T2}.

\vspace{2mm}Similar results hold for the rest of cases listed in
Proposition \ref{p3}. More exactly, the following assertions are
valid:

\vspace{2mm}$\bullet$ $Spec\ D=(1,2,3)$ corresponds to derivations
of type T2 for
$$\beta=3,\ \ \rho=1,\ \ \gamma=\lambda=\mu=0,$$

$\bullet$ $Spec\ D=(1,2,4)$ corresponds to derivations of type T2
for
$$\beta=4,\ \ \rho=1,\ \ \gamma=\lambda=\mu=0,$$

$\bullet$ $Spec\ D=(1,2,a)$ corresponds to derivations of type T2
for
$$\beta=a,\ \ \rho=1,\ \ \gamma=\lambda=\mu=0,$$

$\bullet$ $Spec\ D=(1,a,2a)$ corresponds to derivations of type T2
for
$$\beta=1,\ \ \rho=a,\ \ \gamma=\lambda=\mu=0,$$

$\bullet$ $Spec\ D=(1,a,a+1)$ corresponds to derivations of type
T3 for
$$\beta=1,\ \ \ \zeta=a,\ \ \delta=\eta=0,$$

$\bullet$ $Spec\ D=(1,\lambda,\lambda-1)$ corresponds to
derivations of type T3 for
$$\beta=1,\ \ \ \zeta=\lambda-1,\ \ \delta=\eta=0.$$

Summing the results proved before we get the next classifying
result.

\begin{teo}\label{t1}
Any real 3-dimensional commutative algebra having a semisimple
nonsingular derivation has necessarily a derivation with one of
the spectra
$$1)\ (1,-1,2),\ \ \ \ 2)\ (1,1,2).$$Accordingly, each such algebra is
isomorphic to one the next four algebras:
$$\begin{array}{llllll} 1^{\circ})\ \ e_{1}^{2}=e_{3}&e_{1}e_{2}=0&e_{1}e_{3}=0&e_{2}^{2}=0&e_{2}e_{3}=e_{1}&e_{3}^{2}=0\\
2^{\circ})\ \ e_{1}^{2}=0&e_{1}e_{2}=0&e_{1}e_{3}=0&e_{2}^{2}=0&e_{2}e_{3}=0&e_{3}^{2}=e_{2}\\
3^{\circ})\ \ e_{1}^{2}=0&e_{1}e_{2}=e_{3}&e_{1}e_{3}=0&e_{2}^{2}=0&e_{2}e_{3}=0&e_{3}^{2}=0\\
4^{\circ})\ \
e_{1}^{2}=e_{3}&e_{1}e_{2}=0&e_{1}e_{3}=0&e_{2}^{2}=e_{3}&e_{2}e_{3}=0&e_{3}^{2}=0.
 \end{array}$$
The algebras $1^{\circ} - 4^{\circ}$ are not mutually isomorphic.
\end{teo}

Then, the homogeneous quadratic differential systems on
$\mathbb{R}^{3}$ having a semisimple nonsingular derivation are
classified, up to a center-affine equivalence, in accordance with
the classification, up to an isomorphism, of real 3-dimensional
commutative algebra having a semisimple nonsingular derivation. By
using Theorem \ref{t1} we get the next result.

\begin{teo}\label{t2}
Any homogeneous quadratic differential system on $\mathbb{R}^{3}$
having a semi-simple nonsingular derivation is affinely-equivalent
to one of the following systems:
$$\begin{array}{llll}1)\ \left\{\begin{array}{l}\dot{x}=2yz\\ \dot{y}=0\\ \dot{z}=x^{2} \end{array}\right.&
2)\ \left\{\begin{array}{l}\dot{x}=0\\ \dot{y}=0\\
\dot{z}=x^{2}\end{array}\right. & 3)\ \left\{\begin{array}{l}\dot{x}=0\\
\dot{y}=0\\ \dot{z}=2xy
\end{array}\right.&
4)\ \left\{\begin{array}{l}\dot{x}=0\\ \dot{y}=0\\
\dot{z}=x^{2}+y^{2}. \end{array}\right. \end{array}$$
\end{teo}

\section{Conclusions}

The existence of a semisimple nonsingular derivation for a real
3-dimensional commutative algebra implies the existence of a large
family of symmetries reflected by the existence of suitable bases
where the most part of structure constants become zero (see
Theorem \ref{t1}). In these bases, the corresponding homogeneous
quadratic differential systems get the simplest form (see Theorem
\ref{t2}). This implies the existence of geometric symmetries for
the set of all its integral curves. From algebraic point of view
it follows that the existence of a semisimple nonsingular
derivation involves the solvability of these algebras even the
nilpotence property for most part of them. Moreover, algebras
$A_{2}-A_{4}$ are necessarily associative, while algebra $A_{1}$
is neither associative nor power-associative. Further, each ideal
of algebra $A$ allows to identify a decoupled subsystem of the
associated HQDS. Consequently, from geometrical point of view, let
us remark that all its nonsingular integral curves have zero
torsion, because each such algebra is solvable and there exists at
least a linear prime integral of the associated system.  Moreover,
Theorem \ref{t2} assures that integral curves of systems 2)-4)
have both curvature and torsion tensors zero, so that these curves
are lying on lines parallel to one of the coordinate axes.


\begin{thebibliography}{99}

\bibitem{Bur1} {\sc Burdujan I.}, Quadratic differential systems,
(Romanian), Pim, Ia\c si, 2008.

\bibitem{Bur2} {\sc Burdujan I.}, Homogeneous quadratic dynamical
systems on $\mathbb{R}^{^{3}}$ having derivations with complex
eigenvalues, Libertas Mathematica XXVIII (2008) 69--92.

\bibitem{Bur3} {\sc Burdujan I.}, A classification of a class of
homogeneous quadratic dynamical systems on $\mathbb{R}^{^{3}}$
with derivations, Bull. I.P. Ia\c si,
 Sect. Matematica, Mecanica teoretica, Fizica, LIV (2008) 37--47.

\bibitem{Bur4} {\sc Burdujan I.}, Classification of
Quadratic Differential Systems on $\mathbb{R}^{3}$ having a
nilpotent of order 3 Derivation, Libertas Mathematica  XXIX (2009)
47-64.

\bibitem{Da1}{\sc T. Date}, Classifications and Analysis of Two-dimensional Real
Homogeneous Quadratic Differential Equation Systems, J. Diff. Eqs.
32 (1979) 311--334.

\bibitem{KaYo1} {\sc Kaplan J. L. , Yorke J. A.}, Nonassociative real algebras and
quadratic differential equations, Nonlinear Analysis, Theory,
Methods and Applications, v. \textbf{3},(1), 1979 , 49--51.

\bibitem{Kapl} {\sc I. Kaplansky}, Algebras with many derivations, "Aspects of Mathematics and its
Applications" (ed. J.A. Barroso), Elsevier Science Publishers
B.V., 1986, 431.

\bibitem{KiSa1}{\sc K.M. Kinyon , A.A. Sagle}, Quadratic Dynamical Systems and
Algebras, J. of Diff. Eqs. \textbf{117} (1995) 67--127 .

\bibitem{KiSa2} {\sc Kinyon K. M., Sagle A. A.}, Automorphisms and Derivations of
Differential Equations and Algebras, Rocky Mountain J. of Math.,
\textbf{v. 24}, no. 1, 1994, 135--153.

\bibitem{KKOsh} {\sc Krasnov Y., Kononovich A., Osharovich G.}, On a structure of the fixed point set of homogeneous maps,
 Discrete and Continuous Dynamical Systems Ser.S, v.6, no.4, \textbf{2013},
 1017-1027.

\bibitem{Lang} {\sc Lang S.} Algebra (3rd edition),
Springer-Verlag.
2002.


\bibitem{Ma1} {\sc L. Markus}, Quadratic Differential Equations and Non-associative
Algebras, in ''Contributions to the Theory of Nonlinear
Oscillations'' Annals of Mathematics Studies, no. 45, Princeton
University Press, Princeton, N. Y., 1960.


\bibitem{SW} {\sc A.A. Sagle, R. Walde}, Introduction to Lie groups and Lie
algebras, Academic Press, New York, 1973.

\bibitem{Scha1} {\sc  R.D. Schafer}, An Introduction to Nonassociative
Algebras, Academic Press, New York, London, 1966.

\bibitem{VuSi} {\sc N.I. Vulpe, K.S. Sibirski\u{\i}}, Geometrical
Classification of Qua\-dra\-tic differential systems (Russian),
Differentialnye Uravnenje 13 no. 5 (1977) 803--814.

\bibitem{Wa1} {\sc S. Walker}, Algebras and differential equations, Hadronic
Press, Palm Harbor, 1991.

\end{thebibliography}
 \end{document}